\documentclass[12pt]{article}

\usepackage[leqno]{amsmath}
\usepackage{amsthm,amscd}
\usepackage{amsfonts} 
\usepackage{amssymb} 

\makeatletter
\def\@cite#1#2{{#1\if@tempswa , #2\fi}}
\def\@lbibitem[#1]#2{\item\if@filesw 
      { \def\protect##1{\string ##1\space}\immediate
        \write\@auxout{\string\bibcite{#2}{#1}}}\fi\ignorespaces}
\makeatother
\newtheorem{theorem}{Theorem}[section]
\newtheorem{prop}{Proposition}[section]
\newtheorem{co}{Corollary}[section]
\newtheorem{lm}{Lemma}[section]
\newtheorem{rem}{Remark}[section]
\newtheorem{example}{Example}[section]
\newtheorem{definition}{Definition}[section]
\newtheorem{assumption}{Assumption}[section]
\def\Proof{\noindent{\sl Proof.}\qquad}
\def\QED{\hfill\hbox{\vrule width 4pt height 6pt depth 1.5pt}\par\bigskip}

\def\diag{\mathop{\rm diag}}

\newcommand{\trans}{{^t}}

\def\kmms{\kern-\mathsurround}
\newcommand{\vx}{\kmms\mbox{\boldmath $x$\kmms}}
\newcommand{\vn}{\kmms\mbox{\boldmath $n$\kmms}}
\newcommand{\vz}{\kmms\mbox{\boldmath $z$\kmms}}

\newcommand{\calX}{{\cal X}}
\newcommand{\calY}{{\cal Y}}
\newcommand{\calZ}{{\cal Z}}
\newcommand{\calG}{{\cal G}}
\newcommand{\calN}{{\cal N}}
\newcommand{\calH}{{\cal H}}
\newcommand{\calV}{{\cal V}}
\newcommand{\calF}{{\cal F}}
\newcommand{\calU}{{\cal U}}
\newcommand{\calK}{{\cal K}}

\newcommand{\calL}{{\cal L}}

\newcommand{\bbR}{{\mathbb R}}
\newcommand{\bbS}{{\mathbb S}}

\title{Hierarchical orbital decompositions 
and extended decomposable distributions
}
\author{
Hidehiko Kamiya \\ 
{\it Okayama University} \\ 
and \\
Akimichi Takemura \\
{\it University of Tokyo} 
}
\date{June 2006}

\begin{document}
\maketitle

\begin{abstract}
\noindent
Elliptically contoured distributions can be considered to be 
the distributions for which the contours of the density functions 
are proportional ellipsoids. 
Kamiya, Takemura and Kuriki (2006) generalized the elliptically contoured 
distributions to star-shaped distributions, for which the contours are allowed to 
be arbitrary proportional star-shaped sets. 
This was achieved by considering the so-called orbital decomposition 
of the sample space in the general framework of group invariance. 
In the present paper, we extend their results by conducting the 
orbital decompositions in steps and obtaining a further, hierarchical 
decomposition of the sample space. 
This allows us to construct probability models and distributions with
further independence structures.
The general results are applied to the star-shaped distributions with a 
certain symmetric structure, the distributions related to the two-sample 
Wishart problem and the distributions of preference rankings.

\end{abstract}

\section{Introduction}





Elliptically contoured distributions are defined to be the distributions for which 
the contours of the density functions are proportional ellipsoids. 
As a natural generalization of the multivariate normal distributions,  
they are widely used as a distributional assumption 
(Kelker (1970), Cambanis, Huang and Simons (1981), Fang and Anderson (1990)).  
Fang and Zhang (1990) discuss ``generalized multivariate analysis'' based on 
elliptically contoured distribuitons. 
In the meantime, from another perspective, elliptically contoured distributions can be 
obtained from spherical distributions by affine transformation. 
Extending the $l_2$\kmms-norm in spherical distributions 
to the $l_q$\kmms-norm, 
$q>0,$ Osiewalski and Steel (1993) 
introduced $l_q$\kmms-spherical distributions. 

Generalizing 
elliptically contoured distributions 
and $l_q$\kmms-spherical distributions, 
Kamiya, Takemura and Kuriki (2006) defined the so-called star-shaped 
distributions, for which the contours of the density functions are 
proportional to (the boundaries of) arbitrary star-shaped sets $\tilde{\calZ},$ 
called cross sections (see also Fern\'{a}ndez, Osiewalski and Steel (1995) and 
Ferreira and Steel (2005)).
They 
showed that basic facts about the independence of 
the ``length'' and ``direction'' continue to hold for star-shaped distributions. 
However, when the star-shaped set has a symmetric structure, 
we can make a more detailed investigation into this distribution.

In the star-shaped distribution,  
the cross section $\tilde{\calZ}\subset \bbR^p-\{ {\bf 0} \}$ is allowed to be 
an arbitrary star-shaped set---a set which 
intersects each ray emanating from the origin exactly once, 
and the density is assumed to be constant on each 
proportional star-shaped set 
$g\tilde{\calZ}=\{ g\tilde{\vz}: \tilde{\vz} \in \tilde{\calZ} \}, \ g>0.$  
However, there are some cases where 
we have some symmetry; 
as in the case of elliptically contoured distributions, 
we might be able to assume that 
$\tilde{\calZ}$ is symmetric about the origin: 
$\tilde{\calZ}=-\tilde{\calZ}=\{-\tilde{\vz}: \tilde{\vz} \in \tilde{\calZ} \}.$ 
In those cases, $\tilde{\calZ}$ may be obtained as 
$\tilde{\calZ}=\{ \pm 1\}\calZ=\{ \pm \vz: \vz \in \calZ \}
=\calZ \cup (-\calZ)$ in terms of 
a set $\calZ \subset \bbR^{p}-\{ {\bf 0} \}$ which  
intersects each line through the origin exactly once.
As long as this condition is satisfied, $\calZ$ is allowed to be an 
arbitrary set. 
Now, suppose $\vx \in \bbR^p$ is distributed according to a star-shaped 
distribution with respect to such a symmetric $\tilde{\calZ}=\calZ \cup (-\calZ).$ 
Then the density of this distribution is constant on each 
$g\tilde{\calZ}, \ g>0,$ 
and the distribution of $-\vx$ is the same 
as that of $\vx.$

In the above situation, we cannot deal with the skewness of the 
distributions. 
However, we can go further and consider 
those distributions 
whose  
densities are constant on each $g\calZ, \ g \ne 0,$ but 
not necessarily constant on each 
$g\tilde{\calZ}=g\calZ \cup (-g\calZ), \ g>0,$ 
where $\tilde{\calZ}=\calZ \cup (-\calZ).$ 
That is, 
the value of the density on $g\calZ$ can differ from the value 
on $-g\calZ$ for $g>0.$    
In such a case, distributions of $\vx$ and $-\vx$ are not the
same.

These types of distributions can be studied by decomposing 
$\vx \ (\ne {\bf 0})$ uniquely 
as $\vx=\epsilon h \vz, \ \ \epsilon=\pm 1, \ h>0, \ \vz \in \calZ,$ 
with respect to $\calZ.$ 
As with Kamiya, Takemura and Kuriki (2006), these problems can be 
treated as special cases of a general discussion in terms of abstract group 
invariance, and we choose to do so in this paper. 
This approach enables us to apply the obtained general results to 
the distributions of random matrices and moreover discrete distributions. 

In the general framework of group invariance, the present problem
corresponds to the decomposition of the sample space under the action
of an invariance group, called the orbital decomposition (Wijsman
(1990)).  Conducting this decomposition twice, we obtain a
hierarchical decomposition into three parts.  
When a group $\calG$ acts on
a sample space $\calX,$ any subgroup $\calH$ of $\calG$ acts 
on each $\calG$\kmms-orbit.   
By choosing an appropriate $\calH$ and performing a hierarchical orbital
decomposition, we can construct probability models with the
corresponding hierarchical independence structures.  
This allows us to propose in our general framework 
new probability models 
for various statistical problems, 
as illustrated with modeling of preference rankings in Section 4.3.

%
%



The organization of this paper is as follows. 
In Section 2, we summarize some fundamental facts about  
group actions and orbital decompositions, and review the results 
about the decomposable distributions studied in Kamiya, Takemura 
and Kuriki (2006). 
In Section 3 we introduce a further, hierarchical decomposition 
by means of a subgroup action and define extended decomposable
distributions.  We establish various facts on hierarchical orbital
decompositions and derive distributional properties of extended
decomposable distributions.
In the final section, we apply the general results 
to the star-shaped distributions (Section 4.1), 
the two-sample Wishart problem (Section 4.2) 
and the distributions of preference rankings (Section 4.3). 


\section{Orbital decomposition and decomposable distribution}

In this section we summarize some fundamental facts about  
group actions and orbital decompositions, and review the results 
about the decomposable distributions.
For group invariance in statistics, the reader is referred to 
Eaton (1989), Barndorff-Nielsen, Bl{\ae}sild and Eriksen (1989) 
and Wijsman (1990). 
For global cross sections and orbital decompositions in particular, 
see Wijsman (1967, 1986), Koehn (1970), Bondar (1976) and Kamiya (1996). 

\subsection{Orbital decomposition}

Let a group $\calG$ act on a space $\calX$ 
(typically the sample space)   
from the left 
$(g, x) \mapsto gx : \ \calG \times \calX \to \calX.$ 
We write the action of $\calG$ on $\calX$ as $(\calG, \calX).$ 
 
Let $\calG x = \{gx: g \in \calG \}$ be 
the {\it orbit} containing $x \in \calX,$ and let 
$\calX/\calG = \{ \calG x : x \in \cal X \}$ be the 
{\it orbit space}, i.e., the set of all orbits. 
When $\calX$ consists of a single orbit $\calX = \calG x,$ 
the action is said to be {\it transitive}.

Indicate by 
$\calG_x = \{ g \in \calG : gx = x \}$ 
the {\it isotropy subgroup} at $x \in \calX.$
When $\calG_x = \{ e \}$ for all $x\in \calX,$ 
the action is said to be {\it free,} 
where $e$ denotes the identity element of $\calG.$
In general, the isotropy subgroups at two points on a common orbit 
are conjugate to each other:
$\calG_{gx} = g \calG_x g^{-1}, \ \ g \in \calG, \ x \in \calX.$ 

The set of {\it left cosets} 
$g \calG_x = \{ gg': g'\in \calG_x \}, \ g \in \calG,$ is 
called the {\it left coset space} of $\calG$ modulo $\calG_x,$ 
and is denoted by 
$\calG/\calG_x = \{ g \calG_x: g \in \calG \}.$ 
The {\it canonical map} $\pi: \calG \to \calG/\calG_x$ is defined 
by $\pi(g) = g \calG_x, \ g \in \calG.$
The group $\calG$ or more generally its subgroup $\calH<\calG$ 
acts on $\calG/\calG_x$ by 
\begin{equation}
\label{eq:hgGx}
(h, g\calG_x) \mapsto (hg)\calG_x, \quad h \in \calH, \ g \in \calG.
\end{equation}
This action is not transitive unless $\calH$ includes a complete 
set of representatives of $g\calG_x, \ g \in \calG,$ i.e., 
$\calG=\bigcup_{h\in \calH}h\calG_x.$  

A subset $\calZ \subset \calX$ is said to be a {\it cross section} if 
$\calZ$ intersects each orbit $\calG x, \ x \in \calX,$ 
exactly once. 
So any cross section $\calZ$ is in one-to-one correspondence with 
the orbit space $\calX/\calG$ by 
$z \leftrightarrow \calG z, \ z \in \calZ.$
A cross section $\calZ$ having a common  
isotropy subgroup is called a {\it global cross section}: 
$\calG_z = \calG_0,$ say, for all $z \in \calZ.$ 
Unlike a mere cross section, a global cross section does not always exist. 
A global cross section exists if and only if the isotropy subgroups 
$\calG_x, \ x \in \calX,$ are all conjugate to one another. 

Unless otherwise stated, however, we assume from now on that
there does exist a global cross section $\calZ.$ 
Then, we have 
the following one-to-one correspondence: 
\begin{eqnarray}
\label{eq:X<->YxZ}
\calX & \leftrightarrow & \calY \times \calZ, \\
x & \leftrightarrow & (y, z), \ \ \ x = gz, \ \ y = \pi(g), 
\ \ g \in \calG, \nonumber
\end{eqnarray}
where $\calY=\calG/\calG_0$ with $\calG_0=\calG_z, \ z \in \calZ.$ 
Decomposition \eqref{eq:X<->YxZ} is called the {\it orbital decomposition} 
of $\calX$ (or $x$\kmms) with respect to 
$\calZ.$
In \eqref{eq:X<->YxZ} we can think of $y$ and $z$ as 
functions $y=y(x)$ and $z=z(x)$ of $x.$
Under the action of $\calG$ on $\calX,$ 
$y(x)$ is equivariant and $z(x)$ is invariant:
$y(gx)=gy(x), \ z(gx)=z(x), \ \ g \in \calG, \ x \in \calX.$


\

We move on to reviewing some properties of global cross 
sections obtained in Kamiya, Takemura and Kuriki (2006). 

Let $\calZ$ be a global cross section. 
Then $g \calZ = \{ gz: z \in \calZ \}$ for each $g \in \calG$
is again a global cross section.  
We say $g \calZ, \ g \in \calG,$ are {\it proportional to} $\calZ,$ and 
call $\{ g \calZ : g \in \calG \}$ the {\it family of proportional
global cross sections}.
In $\calX=\bigcup_{g\in \calG}g\calZ,$ it holds that   
$g_1 \calZ \cap g_2 \calZ \ne \emptyset$ for $g_1, g_2 \in \calG$  
implies $g_1 \calZ = g_2 \calZ,$ so the family of proportional 
global cross sections gives a partition of $\calX.$

From a given global cross section $\calZ,$ we can construct 
a general cross section $\calZ'$ by 
changing 
the points of $\calZ$ within their orbits.  
For $\calZ'$ to be global, i.e., for the
isotropy subgroups to be the same on the whole of $\calZ',$ 
these 
changes 
of the points have to be made 
subject to some restriction as follows.  
Let
$\calN=\{ g \in \calG: g \calG_0 g^{-1}=\calG_0 \}$
be the {\it normalizer} of the common isotropy subgroup $\calG_0$ 
on $\calZ.$ 
Then a subset $\calZ' \subset \calX$ is a global cross section if and only if 
it can be written as 
\begin{equation}
\label{eq:g0nz}
\calZ'=\{ g_0 n_z z: z \in \calZ \}
\end{equation}
for some $g_0 \in \calG$ and $n_z \in \calN, \ z \in \calZ.$ 
%

Under the change from $\calZ$ to $\calZ'$ in \eqref{eq:g0nz}, 
the equivariant part transforms as follows. 
Let $x \leftrightarrow (y,z)$ be the orbital decomposition
with respect to $\calZ,$ 
and 
let $ x \leftrightarrow (y', z')$ be the orbital decomposition
with respect to 
the $\calZ'$ in (\ref{eq:g0nz}). 
Then we have 
\begin{equation}
\label{eq:transformation of y}
y' = y n_z^{-1} g_0^{-1}.
\end{equation}
%

\subsection{Decomposable distribution}
\label{sec:decomposable-distribution}

In this subsection, we review the decomposable
distributions defined in Kamiya, Takemura and Kuriki (2006). 


Throughout the rest of the paper, we make the following assumptions: 
(a) $\calX$ is a locally compact Hausdorff space; 
(b) $\calG$ is a second countable, locally compact Hausdorff 
topological group acting continuously on $\calX;$ 
(c) $\calG_0$ is compact; 
and 
(d) $\calZ$ is locally compact and the bijection
$x \leftrightarrow (y, z)$ with respect to $\calZ$ 
is bimeasurable.

We consider distributions on $\calX$ which have densities $f(x)$ 
with respect to a dominating measure $\lambda.$ 
Measure $\lambda$ is assumed to be  
relatively invariant with multiplier $\chi:$ 
$\lambda(d(gx))=\chi(g)\lambda(dx), 
\ g \in \calG.$
%
Then, we say 
a distribution $f(x)\lambda(dx)$ is a {\it decomposable distribution} 
with respect to $\calZ$ 
if it is of the form 
$f(x)\lambda(dx) = f_{ \calY }(y(x))f_{ \calZ }(z(x))\lambda(dx).$ 
In particular, we say it is {\it cross-sectionally contoured} if 
$f_{ \calZ }(z)$ is constant, and {\it orbitally contoured} if 
$f_{ \calY }(y)$ is constant. 
We mainly study cross-sectionally contoured distributions 
because a decomposable distribution 
$f_{ \calY }(y(x)) f_{ \calZ }(z(x)) \lambda(dx)$
can always be thought of as a cross-sectionally contoured distribution 
with density $f_{ \calY }(y(x))$ with respect to 
$\tilde{\lambda}(dx) := f_{ \calZ }(z(x))\lambda(dx)$
Obviously, a distribution $f(x)\lambda(dx)$ is
cross-sectionally contoured with respect to $\calZ$ 
if and only if $f(x)$ is constant on each
proportional global cross section $g\calZ, \ g \in \calG.$


Topological assumption (b) about $\calG$ implies that 
there exists a left Haar measure $\mu_{\calG}$ on $\calG,$ 
which is unique up to a multiplicative constant.  
By the compactness of $\calG_0$ assumed in (c), 
we have the induced measure 
$\mu_{ \calY } = \pi(\mu_{ \calG }) = \mu_{ \calG } \pi^{-1}$ 
on $\calY$ (Proposition 2.3.5 and Corollary 7.4.4 of Wijsman (1990)).
Again by the same assumption (c), we can define 
$\bar{\chi}(y), \ y \in \calY,$ by 
$\bar{\chi}(y) = \chi(g)$ with $g \in \pi^{-1}(\{ y \}).$
By abuse of notation, 
we will write $\chi(y)$ for $\bar{\chi}(y).$

In terms of these, 
$\lambda(dx)$ is factored as   
\begin{equation}
\label{eq:decomposition of lambda}
\lambda(dx) = \chi(y)\mu_{ \calY }(dy) \nu_{ \calZ }(dz)
\end{equation}
(Theorem 7.5.1 of Wijsman (1990), 
Theorem 10.1.2 of Farrell (1985)). 
Here, we are identifying $\calX$ with $\calY \times \calZ.$
Existence of a density $f(x)$ with respect to $\lambda$ implies 
$\nu_{\calZ}$ is a finite measure, so from now on 
we assume that $\nu_{ \calZ }(dz)$ is 
standardized to be a probability measure on $\calZ.$ 
From \eqref{eq:decomposition of lambda} 
we immediately obtain the following result. 

\begin{prop} {\bf (Kamiya, Takemura and Kuriki (2006))} \quad 
\label{prop:distributions of y and z}
Suppose that $x$ is distributed according to a cross-sectionally
contoured distribution $f_{\calY}(y(x))\lambda(dx).$
Then we have:
\begin{enumerate}
\item $y = y(x)$ and $z = z(x)$ are independently distributed.
\item The distribution of $y$ is 
$f_{\calY}(y)\chi(y)\mu_{ \calY }(dy).$ 
\item The distribution of $z$ does not depend on $f_{ \calY }.$ 
\end{enumerate}
\end{prop}

Note that since $\nu_{\calZ}(dz)$ is taken to be a probability measure, 
we do not need a normalizing constant in $f_{\calY}(y)\chi(y)\mu_{ \calY }(dy).$  
(To put it another way, the version of $\mu_{\calY}$ is taken in this way.) 

\section{Hierarchical orbital decomposition and extended decomposable 
distribution}

In this section we introduce a further, hierarchical decomposition 
and define extended decomposable distributions.

\subsection{Hierarchical orbital decomposition}


In this subsection, we give a further factorization of the $\calG$\kmms-orbital
decomposition. 
This is obtained by decomposing the equivariant part $\calG/\calG_0$ 
by means of the action of a subgroup 
of $\calG.$ 

We continue to assume that there exists a global cross section $\calZ$ 
with the common isotropy subgroup $\calG_0.$ 
Furthermore, let $\calH$ be a subgroup of $\calG.$ 

As in Section 2.1, we have the decomposition
\begin{equation}
\label{decomposition of X}
\calX \leftrightarrow \calG/\calG_0 \times \calZ.
\end{equation}
Now, $\calH$ acts on $\calG/\calG_0$
by 
\eqref{eq:hgGx} with $\calG_x=\calG_0.$
Note that instead of this action 
we may equivalently consider the action
of $\calH$ on $\calG z_0:
(h, gz_0) \mapsto (hg)z_0, \ z_0 \in \calZ.$ 
In particular, we have $\calH_{g\calG_0}=\calH_{gz_0}, \ g \in \calG,$  
from which we obtain the following lemma: 
\begin{lm} \quad 
\label{lm:H,free}
Action $(\calH, \calX)$ is free if and only if action 
$(\calH, \calG/\calG_0)$ is free. 
\end{lm}

Now suppose a global cross 
section $\calV \subset \calG/\calG_0$ exists for action $(\calH, \calG/\calG_0).$  
The existence of a global cross section $\calV$
leads to a further decomposition of
(\ref{decomposition of X}) as follows.

Denote the common isotropy subgroup at the points of
$\calV$ by $\calH_0.$
Then $\calG/\calG_0$ is decomposed as
\begin{equation}
\label{decomposition of G/G_0}
\calG/\calG_0 \leftrightarrow \calH/\calH_0 \times \calV.
\end{equation}
We can take $\calV$ in such a way that $\calG_0 \in \calV;$
in that case, we can write $\calH_0$ as
\[
\calH_0 = \calH_{ \calG_0 } = \{ h \in \calH: h \calG_0 = \calG_0 \}
= \calH \cap \calG_0.
\]
From now on, we always take $\calV$ in this way.

Combining (\ref{decomposition of X}) and
(\ref{decomposition of G/G_0}), we have the decomposition
\begin{equation}
\label{eq:decomposition-into-3parts}
\calX \leftrightarrow \calH/\calH_0 \times \calV \times \calZ.
\end{equation}
Our questions are:
\begin{enumerate}
\setlength{\itemsep}{0pt}
\item[(i)] specifying the condition for $\calV$ to exist, and
\item[(ii)] expressing $\calV$ in a concrete form.
\end{enumerate}
Note that 
the orbits under $(\calH, \calG/\calG_0)$ are of the form 
\[
\left\{ hg\calG_0: h \in \calH \right\}
\subset \calG/\calG_0, 
\quad g\in \calG.
\]
This suggests that the above questions are 
closely related to the properties of the
double cosets 
\begin{eqnarray*}
\calH g \calG_0
&=&\{ hgg_0: h \in\calH, \ g_0\in \calG_0 \} \\ 
&=&\pi^{-1}\left( 
\left\{ hg\calG_0: h \in \calH \right\} 
\right)\subset \calG, 
\quad g \in \calG, 
\end{eqnarray*}
in $\calG.$
The following lemma indicates this fact.

\begin{lm}
\label{lm:double coset representatives}
Let $\calG'
\subset \calG.$
Then $\calV= \{ g' \calG_0 : g' \in \calG' \}\subset\calG/\calG_0$ 
is a cross section
for the action of $\calH$ on $\calG/\calG_0$ 
such that $g'\calG_0\ne g''\calG_0$ for $g'\ne g'', \ g', g''\in \calG',$ 
if and only if
$\calG'$ 
is a complete set of representatives
of the double cosets $\calH g \calG_0, \ g \in \calG,$ in $\calG:$ 
\[
\calG=\bigsqcup_{g'\in \calG'}\calH g' \calG_0 \quad \text{{\rm (disjoint union).}} 
\]
\end{lm}
\Proof
{\it Necessity}:   
Suppose $\calV= \{ g' \calG_0 : g' \in \calG' \}$ is a cross section 
for action $(\calH, \calG/\calG_0)$ and that 
$g'\calG_0 = g''\calG_0$ for $g', g''\in \calG'$ implies $g'=g''.$ 
We want to prove $\calG=\bigsqcup_{g'\in \calG'}\calH g' \calG_0.$ 
It suffices to verify (a)$\calG\subset \bigcup_{g'\in \calG'}\calH g' \calG_0;$ 
and (b)$\calH g' \calG_0 = \calH g'' \calG_0$ for $g', g''\in \calG'$ implies $g'=g''.$ 
First, (a) is shown as follows. 
Let $g$ be an arbitrary element of $\calG.$ 
Then, since $\calV$ intersects the orbit containing $g\calG_0\in\calG/\calG_0$ 
at least once, there exist $h\in \calH$ and $g'\in \calG'$ such that 
$g\calG_0=hg'\calG_0.$ 
Thus $g$ can be written as $g=hg'g_0$ with some $g_0\in \calG_0.$  
Therefore, $g\in\calH g' \calG_0 \subset \bigcup_{g''\in \calG'}\calH g'' \calG_0.$ 
Next, (b) is proved as follows. 
Suppose $\calH g' \calG_0 = \calH g'' \calG_0$ for $g', g''\in \calG'.$ 
Then $g'=hg''g_0$ for some $h\in \calH$ and $g_0\in \calG_0,$ and thus 
we have $g'\calG_0=hg''\calG_0.$ 
Now, since $\calV$ intersects the orbit containing $g''\calG_0$ at most once, 
we obtain $g'\calG_0=g''\calG_0.$ 
Therefore, we get $g'=g''$ by our assumption. 

{\it Sufficiency}:   
Suppose $\calG=\bigsqcup_{g'\in \calG'}\calH g' \calG_0.$ 
We want to show (a)$\calV= \{ g' \calG_0 : g' \in \calG' \}$ intersects 
each orbit under $(\calH, \calG/\calG_0)$ at least once; 
(b)$\calV$ intersects each orbit at most once; and 
(c)$g'\calG_0 = g''\calG_0$ for $g', g''\in \calG'$ implies $g'=g''.$  
We begin by showing (a). 
Let $g \calG_0 \in \calG/\calG_0, \ g \in \calG,$ be arbitrarily given. 
Pick any $g_1\in g\calG_0\subset
\calG=\bigcup_{g'\in \calG'}\calH g' \calG_0.$ 
Then $g_1$ can be written as $g_1=hg'g_0$ for some $h\in \calH, \ 
g'\in\calG'$ and $g_0\in \calG_0.$ 
Hence $g\calG_0=g_1\calG_0=hg'\calG_0.$ 
This shows that the orbit containing $g\calG_0$ intersects $\calV$ 
at least once. 
Next we verify (b). 
Suppose $g'\calG_0=hg''\calG_0$ for $g', g''\in \calG'$ and $h\in \calH.$ 
Then $\calH g' \calG_0=\calH g'' \calG_0,$ which implies $g'=g''$ since 
$\calG=\bigsqcup_{g'\in \calG'}\calH g' \calG_0$ is a disjoint union. 
Hence we have $g'\calG_0=g''\calG_0.$ 
This observation shows (b). 
Finally, (c) can be verified similarly as (b).   
\QED


Now, concerning the existence of $\calV,$ we state the following
theorem.

\begin{theorem} \quad
\label{th:existence of V}
Suppose that there exists a global cross section $\calZ$ 
for the action of $\calG$ on $\calX,$ with the common 
isotropy subgroup denoted by $\calG_0.$ 
Let $\calG'$ be a complete set of
representatives of 
the double cosets $\calH g \calG_0, \ g \in \calG,$ 
in $\calG.$ 
Then a global cross section $\calV$ exists for the action 
of $\calH$ on $\calG/\calG_0$ 
if and only if 
\[
  \calH \cap g' \calG_0 g'^{-1}, \quad g' \in \calG', 
\]
are all conjugate to one another 
in $\calH.$
\end{theorem}

\Proof
First note that a global cross section $\calV$ exists for action
$(\calH, \calG/\calG_0)$ if and only if the isotropy subgroups
$\calH_{g \calG_0}, \ g \in \calG,$ are all conjugate in $\calH.$

Every 
$g \in \calG = \bigsqcup_{ g' \in \calG' } \calH g' \calG_0$ 
can be written 
in the form 
$g = h g' g_0, \ h \in \calH, \ g' \in \calG', \ g_0 \in \calG_0,$ 
and thus we have 
$\calH_{g \calG_0} = \calH_{h g' g_0 \calG_0} 
=h \calH_{g' \calG_0}h^{-1}.$
Therefore, $\calV$ exists if and only if
$\calH_{g' \calG_0}, \ g' \in \calG',$ are all conjugate in $\calH.$
Here we can write $\calH_{g' \calG_0}$ as
$\calH_{g' \calG_0}
= \{ h \in \calH: h g' \calG_0 = g' \calG_0 \} 
= \{ h \in \calH: g'^{-1}hg' \in \calG_0 \} 
= \calH \cap g' \calG_0 g'^{-1}.$
\QED

\begin{rem} \quad
\label{rem:representatives}
The condition that all $\calH \cap g' \calG_0 g'^{-1}, \ g' \in \calG',$
be conjugate does not depend on the choice of a complete set $\calG'.$
\end{rem}

%


Let us now move on to the second problem---expressing $\calV$
in a 
concrete form. 
The following theorem gives a useful  explicit expression of $\calV.$ 
We omit the proof because it is a simple consequence of 
Lemma \ref{lm:double coset representatives}.

\begin{theorem} \quad
\label{th:form of V}
Suppose that there exists a global cross section $\calZ$ for 
the action of $\calG$ on $\calX,$ with the common isotropy 
subgroup $\calG_0.$
Suppose furthermore that there exists a complete set 
$\calG' = \{ g_i: i\in I \}
\subset \calG$ 
of representatives of the double cosets 
$\calH g\calG_0, \ g \in \calG,$ in $\calG$ 
such that 
$\calH_{g_i \calG_0 }$ 
does not depend on $i \in I.$
Then
\[
\calV = \{ g_i \calG_0: i \in I \}
\]
is a global cross section for the action of 
$\calH$ on $\calG/\calG_0.$
\end{theorem}

\begin{rem} \quad
\label{rem:condition-of-G'}
When the action of $\calH$ on $\calG/\calG_0$ is free, 
any complete set 
of representatives of 
$\calH g \calG_0, \ g \in \calG,$ in $\calG$ satisfies the condition 
of $\calG'$ in Theorem \ref{th:form of V}: 
$\calH_{g_i\calG_0}=\{ e \}$ for all $i\in I.$ 
\end{rem}


Under the assumption of Theorem \ref{th:form of V},
$\calX$ is decomposed as 
\eqref{eq:decomposition-into-3parts}. 
We now prove that $\calV \times \calZ$ 
in \eqref{eq:decomposition-into-3parts} 
is a global cross section for action $(\calH, \calX).$

\begin{theorem} \quad
\label{th:V times Z}
Suppose that there exists a global cross section $\calZ$ for 
the action of $\calG$ on $\calX,$ with the common isotropy 
subgroup $\calG_0.$
Suppose moreover that there exists a complete set 
$\calG'=\{ g_i: i\in I \}$ of representatives 
of the double cosets $\calH g\calG_0, \ g \in \calG,$ in $\calG$ satisfying
the condition of Theorem \ref{th:form of V}, and let
$\calV = \{ g_i \calG_0: i \in I \}.$
Then $\calV \times \calZ$ is in one-to-one correspondence with 
\begin{equation}
\label{eq:decomposable global cross section}
\tilde{\calZ}:=
\calG' \calZ = \{ g_i z: \ i \in I, \ z \in \calZ \},
\end{equation}
and $\tilde{\calZ}$ is a global cross section for the action of $\calH$
on $\calX.$
\end{theorem}

\Proof
It is easy to see that the correspondence
$(g_i \calG_0, z) \leftrightarrow g_iz$ 
between $\calV \times \calZ$ and $\calG' \calZ$ is a 
bijection.
We show below that $\calG' \calZ$ is a global cross section for
action $(\calH, \calX).$

Let $x$ be an arbitrary element of $\calX.$ 
Then $x$ can be written as $x=gz, \ g \in \calG, \ z \in \calZ.$
Furthermore, this $g$ 
can be written as $g=hg_ig_0$ for some 
$h \in \calH, \ i \in I$ and $g_0 \in \calG_0.$  
Hence, $x=hg_ig_0 z=hg_iz$ and so $\calH x = \calH g_i z \ni g_iz.$  
This implies that 
$\calG' \calZ$ intersects each 
orbit $\calH x, \ x \in \calX,$ under $(\calH, \calX)$ at least once.

Next we show that $\calG' \calZ$ intersects each
$\calH x, \ x \in \calX,$ at most once.
Suppose that
$
  hg_iz \in \calG'\calZ
$
for $h \in \calH, \ i \in I$ and $z \in \calZ.$ 
Then there exist 
$i'\in I$ 
and
$z' \in \calZ$ such that $hg_i z = g_{i'}z'.$
Since $\calZ$ is a cross section for action $(\calG, \calX),$
we have $z=z'$ and $hg_i \calG_0 = g_{i'} \calG_0,$ 
which implies $\calH g_i\calG_0=\calH g_{i'}\calG_0$ and thus 
$i=i'.$  
Therefore, we obtain
$
  hg_i z = g_i z.
$
This observation shows that $\calG' \calZ$ intersects each
$\calH x, \ x \in \calX,$ at most once.

It remains to be proved that the isotropy subgroups
$\calH_{g_i z}$ at the points $g_i z \in \calG' \calZ$
are all common. 
But this is obvious because 
$\calH_{g_i z} = \calH_{g_i \calG_0}$ 
does not depend on $i\in I$ by the assumption of 
Theorem \ref{th:form of V}. 
\QED

We call a global cross section $\tilde{\calZ}$ 
for action $(\calH, \calX)$ 
of the form
\eqref{eq:decomposable global cross section} 
a {\it decomposable global cross section}. 
Of course, a general global cross section 
for $(\calH, \calX)$ is not necessarily decomposable. 


\

Let us delve into Theorems \ref{th:form of V} and \ref{th:V times Z} in 
two specific cases.

First, consider the case where the action of $\calG$ 
on $\calX$ is free.
In that case, we want to decompose $\calG$ by considering 
the action of $\calH$ on $\calG: (h, g) \mapsto hg, \ 
h \in \calH, \ g \in \calG.$

\begin{co} \quad
\label{co:form of V when (G,X) is free}
The action of $\calH$ on $\calG$ is free, and any complete set
$\{ g_i: i \in I \}$ of representatives of the right cosets 
$\calH g, \ g \in \calG,$ in $\calG$ 
is a cross section for this action.
\end{co}

\Proof
It is trivial to see that action $(\calH, \calG)$ is free. 
The rest is obvious from Theorem \ref{th:form of V} 
with $\calG_0=\{ e \}.$
\QED

The latter statement of the corollary is also apparent  
from the fact that $\calH g, \ g \in \calG,$ are the orbits 
under $(\calH, \calG)$ and 
$\calG=\bigsqcup_{i\in I}\calH g_i.$

In the case of Corollary \ref{co:form of V when (G,X) is free},
$\calG$ is decomposed as
\begin{eqnarray*}
  \calG & \leftrightarrow & \calH \times \{ g_i: i \in I \} \\
    & \leftrightarrow & \calH \times \calH \backslash \calG,
\end{eqnarray*}
where $\calH \backslash \calG$ is the right coset space 
\[
\calH \backslash \calG := \{ \calH g: g \in \calG \}
= \{ \calH g_i: i \in I \}. 
\]

\begin{co} \quad
\label{co:G'Z}
Suppose that the action of $\calG$ on $\calX$ is free,
and let $\calZ$ be a cross section for this action.
Then the action of $\calH$ on $\calX$ is free,
and for any complete set
$\calG' = \{ g_i: i \in I \}$
of representatives of the right cosets $\calH g, \ g \in \calG,$ in $\calG,$ 
the set 
$\tilde{\calZ}=\calG' \calZ = \{ g_i z: i \in I, \ z \in \calZ \}$ is a
cross section for this action.
\end{co}

\Proof
It is evident that the action of $\calH$ on $\calX$ is free. 
The rest follows immediately from 
Theorem \ref{th:V times Z} with 
Lemma \ref{lm:H,free} and 
Remark \ref{rem:condition-of-G'}. 
\QED

As an example, consider the actions related to the star-shaped
distributions---the actions of
$\calG = \bbR^*_\times$ (the multiplicative group of 
nonzero real numbers)
and $\calH = \bbR_+^*$ (the multiplicative group of 
positive real numbers) 
on $\calX = \bbR^p - \{ {\bf 0} \}$ by scalar multiplication. 

In that case, $\calG$ acts on $\calX$ freely.
Moreover, $\{ \pm 1 \}$ is a complete set of representatives
of $\calH g, \ g \in \calG,$ and is thus a cross section for action
$(\calH, \calG)$ by Corollary \ref{co:form of V when (G,X) is free}. 
Accordingly, we have a one-to-one correspondence
\[
  \bbR^*_\times \leftrightarrow \bbR_+^* \times \{ \pm 1 \}.
\]

Furthermore, we have by Corollary \ref{co:G'Z} that
$\tilde{\calZ}=\calG' \calZ=\calZ \cup (-\calZ)$ 
with $\calG' = \{ \pm 1 \}$ is a cross
section for the action of $\calH=\bbR_+^*$ on $\calX=\bbR^p - \{ {\bf 0} \}.$
Let us take $\calZ$ as
\[
  \calZ = \left\{ (x_1, \ldots, x_p)^\trans \in \bbS^{p-1}: x_p > 0 \right\}
      \cup
      \left\{ (x_1,\ldots,x_{p-1}, 0)^\trans \in \bbS^{p-1}:
      (x_1, \ldots, x_{p-1})^\trans \in \calZ_{-1} \right\},
\]
where $\bbS^{p-1}$ denotes the$(p-1)$\kmms-dimensional unit sphere and  
$\calZ_{-1}$ is a cross section for the action of
$\calG=\bbR_{\times}^*$ on $\bbR^{p-1} - \{ {\bf 0} \}, \ {\bf 0} \in \bbR^{p-1}.$
Since it is clear that $\calG' \calZ = \bbS^{p-1}$ is true for $p=2,$
we see by induction on $p$ that for all $p,$
\begin{eqnarray*}
\tilde{\calZ}&=&
\calG' \calZ \\ 
  &=& \left\{ (x_1,\ldots,x_p)^\trans \in \bbS^{p-1}: x_p \ne 0 \right\}
      \cup
      \left\{ (x_1, \ldots, x_{p-1}, 0)^\trans \in \bbS^{p-1} \right\} \\
  &=& \bbS^{p-1},
\end{eqnarray*}
which is clearly a cross section for the action of $\calH=\bbR_+^*$ on
$\calX=\bbR^p - \{ {\bf 0} \}.$
In fact, we can take any $\calZ \subset \calX=\bbR^p-\{ {\bf 0} \}$ which 
intersects each line through the origin in exactly one point. 

\

Next we treat the case which covers the two-sample Wishart problem.

\begin{co} \quad
\label{co:form of V when g=hk}
Let $\calG_0$ be a subgroup of $\calG.$
Suppose that there exists a subgroup $\calK$ of $\calG$ satisfying the
following conditions:
\begin{enumerate}
  \item[{\rm (i)}] Every $g \in \calG$ can be written uniquely in the form
        $g = hk, \ h \in \calH, \ k \in \calK.$
  \item[{\rm (ii)}] $\calG_0$ is a subgroup of $\calK.$
\end{enumerate}
Then, the action of $\calH$ on $\calG/\calG_0$ is free, and
\[
  \calV = \calK/\calG_0 = \{ k \calG_0: k \in \calK \}
\]
is a cross section for this action.
\end{co}

\Proof
Noting that $\calH \cap \calK = \{ e \}$ by assumption (i) and 
that $k \calG_0 k^{-1} \subset \calK, \ k \in \calK,$ by 
assumption (ii),
we have $\calH \cap k \calG_0 k^{-1} = \{ e \}$ for 
any $k \in \calK.$
Therefore, $\calH_{g \calG_0} = \calH \cap g \calG_0 g^{-1}$ 
is trivial for any $g$ in $\calK$ and thus for any $g$ in $\calG:$
\[
  \calH_{g \calG_0} = \calH_{hk \calG_0} = h \calH_{k \calG_0}h^{-1}
  = h\{ e \} h^{-1}= \{ e \}, \ \ \ g=hk, \ h \in \calH, \ k \in \calK.
\]
Hence, action $(\calH, \calG/\calG_0)$ is free.

Let $\calG' \subset \calK$ be a complete set of representatives
of the left cosets $k\calG_0, \ k\in\calK,$ in $\calK: \ 
\calK=\bigsqcup_{g'\in\calG'}g'\calG_0.$ 
Then $\calK/\calG_0=\{ g'\calG_0: g' \in \calG' \},$ and  
by Theorem \ref{th:form of V} 
it suffices to show that $\calG'$ is a
complete set of representatives of the double cosets
$\calH g \calG_0, \ g \in \calG,$ in $\calG: \ 
\calG=\bigsqcup_{g'\in\calG'}\calH g'\calG_0.$ 
Since
$\bigcup_{g'\in\calG'} \calH g' \calG_0 
=\calH (\bigcup_{g'\in\calG'} g' \calG_0) 
= \calH \calK = \calG,$
it remains to show that
\begin{equation}
  \label{eq:g'g''}
  \calH g^{ \prime } \calG_0 = \calH g^{ \prime \prime } \calG_0,
  \ \ \ g^\prime, g^{ \prime \prime } \in \calG^\prime,
\end{equation}
implies $g^\prime = g^{ \prime \prime }.$
Suppose that \eqref{eq:g'g''} holds. 
Then, there exist $h \in \calH$ and $g_0 \in \calG_0$ such that
$g^\prime = h g^{ \prime \prime } g_0.$
By assumptions (i) and (ii), we have
$g^\prime = g^{ \prime \prime } g_0$ 
and thus 
$g^\prime \calG_0 
= g^{ \prime \prime } \calG_0.$
The definition of $\calG'$ implies  
$g^\prime = g^{ \prime \prime }.$
\QED

\begin{rem} \quad 
\label{rem:twoG's}
As was seen in the proof of Corollary \ref{co:form of V when g=hk}, 
if $\calG' \subset \calK$ is a complete set of representatives
of the left cosets $k\calG_0, \ k\in\calK,$ in $\calK,$ then 
$\calG'$ is also a complete set of representatives of the double cosets
$\calH g \calG_0, \ g \in \calG,$ in $\calG.$ 
As can be shown in a similar manner, 
if we assume $\calG'\subset \calK,$  
the converse is true as well 
and the two conditions are in fact equivalent. 
\end{rem}

In the case of Corollary \ref{co:form of V when g=hk},
$\calG/\calG_0$ is decomposed as
\[
  \calG/\calG_0 \leftrightarrow \calH \times \calK/\calG_0.
\]

\begin{co} \quad
\label{co:KZ}
Suppose that there exists a global cross section $\calZ$ for the action
of $\calG$ on $\calX,$ with the common isotropy subgroup denoted by
$\calG_0.$
Suppose moreover that there exists a subgroup $\calK$ satisfying
conditions {\rm (i)} and {\rm (ii)} in Corollary \ref{co:form of V when g=hk}.
Then, the action of $\calH$ on $\calX$ is free, 
and $\tilde{\calZ}=\calK \calZ$ is a cross section for this action.
\end{co}

\Proof 
By Lemma \ref{lm:H,free} and 
Corollary \ref{co:form of V when g=hk}, 
action $(\calH, \calX)$ is free. 
We show below that $\calK \calZ$ is a cross section for $(\calH, \calX).$ 
 
Let $\calG' \subset \calK$ be some complete set of representatives of 
the left cosets $k\calG_0, \ k \in \calK,$ in $\calK.$ 
Then 
by Remark \ref{rem:twoG's}, 
$\calG'$ is a complete set of representatives of 
the double cosets $\calH g \calG_0, \ g \in \calG,$ in $\calG$  
as well. 
Theorem \ref{th:V times Z} implies that $\tilde{\calZ}=\calG'\calZ$ is a 
cross section for $(\calH, \calX).$   
Thus, the proof will be finished if we verify that $\calG'\calZ=\calK\calZ.$ 
But this follows from $\calK=\bigsqcup_{g'\in\calG'}g'\calG_0=\calG'\calG_0.$  
%
%
\QED

As an example, consider the situation related to 
the two-sample Wishart problem. 

Let $\calG=GL(p)$ (the general linear group) 
and
\begin{align}
\label{eq:wishart sample space}
& \calX= \{ (W_1, W_2) \in PD(p) \times PD(p): \\
& \qquad \qquad 
   \text{the $p$ roots of} \ \det(W_1- \lambda (W_1+W_2))=0 \ 
   \text{are all distinct} \}. \nonumber 
\end{align}
The action is
$(B, \ (W_1, W_2)) \mapsto (BW_1B^\trans, \ BW_2B^\trans), \ B \in GL(p).$
Let us take
\begin{equation}
\label{eq:wishart cross section}
  \calZ= \left\{ (\Lambda, I_p- \Lambda): 
  \Lambda = \diag (\lambda_1,\ldots,\lambda_p),
  \ 1 > \lambda_1 > \cdots > \lambda_p > 0 \right\},
\end{equation}
where $I_p$ denotes the $p \times p$ identity matrix.
Then we have
\begin{equation}
\label{eq:G0 for wishart}
  \calG_0 = \left\{ \diag (\epsilon_1,\ldots,\epsilon_p):  
  \epsilon_1 = \pm 1,\ldots,\epsilon_p = \pm 1 \right\}. 
\end{equation}

Now, as a subgroup of $\calG$ consider  $\calH = LT(p),$ 
the group of $p\times p$ lower triangular matrices 
with positive diagonal elements. 
Then the orthogonal group $O(p)$ can serve as the $\calK$ in Corollary
\ref{co:form of V when g=hk}. 
With this $\calK,$ 
\[
\calV=\calK/\calG_0=\{ C\calG_0: C \in O(p) \}
\] 
is the set of 
$p\times p$ orthogonal matrices with the sign of each column ignored. 

Corollary \ref{co:KZ} implies that the action of
$\calH = LT(p)$ on $\calX$ is free and that $\tilde{\calZ}=\calK \calZ$
with $\calK = O(p)$ is a cross section
for $(\calH, \calX).$
We can write $\tilde{\calZ}$ as
\begin{eqnarray*}
\tilde{\calZ}&=&\calK \calZ \\ 
 &=& \left\{ (C \Lambda C^\trans, \ I_p - C \Lambda C^\trans): C \in O(p),
     \ \Lambda = \diag (\lambda_1,\ldots,\lambda_p), \ 
     1 > \lambda_1 > \cdots > \lambda_p > 0 \right\} \\
 &=& \left\{ (U, \ I_p - U): O < U < I_p, \ 
     \text{and the eigenvalues of} \ U \ 
     \text{are all distinct} \right\},
\end{eqnarray*}
where $O$ denotes the null matrix and 
$A<B$ means that $B-A$ is positive definite for symmetric matrices 
$A$ and $B.$ 

\subsection{Extended decomposable distribution}


In this subsection, we discuss the distributional aspect of the 
hierarchical decompositions in the preceding subsection.

Suppose that there exists a subgroup $\calL$ of $\calG$ of the form
\begin{equation}
  \label{eq:calL condition}
  \calL = \calG' \calG_0,
\end{equation}
where $\calG'=\{ g_i: i\in I \}$ 
is a complete set of representatives of the double
cosets $\calH g \calG_0, \ g \in \calG,$ in 
$\calG$ such that $\calH_{g_i \calG_0}$
does not depend on $i \in I.$ 

Then we have 
\[
\calG=\calH \calG' \calG_0=\calH \calL.
\]
Therefore, every $g \in \calG$ can be written in the form
$g = h l \ (\text{or} \ g=hl^{-1}), \ h \in \calH, \ l \in \calL.$
Moreover, by considering the transitive action  
$((h, l), g) \mapsto hgl^{-1}, \ h\in\calH, \ l\in \calL, \ g\in \calG,$ 
of the product group
$\calH \times \calL$ on
$\calG,$ 
we have a bijection
\[
\calG \leftrightarrow (\calH \times \calL)/\calF^*,
\]
where 
$\calF^* = \left\{ (g, g): g \in \calF= \calH \cap \calL \right\}$
is the isotropy subgroup at $e \in \calG.$

\

Before proceeding further, 
let us see the two specific cases considered in the 
preceding subsection.

First, consider the situation in Corollaries
\ref{co:form of V when (G,X) is free} and \ref{co:G'Z}.
Suppose that the action of $\calG$ on $\calX$ is free and that
a complete set $\calG'$ 
of representatives of the right
cosets $\calH g, \ g \in \calG,$ 
forms a subgroup of $\calG.$
Then $\calG'$ can serve as $\calL.$ 
For instance, consider the example immediately after Corollary
\ref{co:G'Z}---the actions related to the star-shaped distributions.
Then $\calG'=\{ \pm 1 \}$ forms a subgroup of $\calG = \bbR^*_\times$ and 
thus can serve as $\calL.$

Next, consider the situation in Corollaries
\ref{co:form of V when g=hk} and \ref{co:KZ}.
Suppose that a subgroup $\calK$ of $\calG$ satisfies 
conditions (i) and (ii) of Corollary \ref{co:form of V when g=hk}, 
and let $\calG'\subset \calK$ be a complete set of representatives of the 
left cosets $k\calG_0, \ k \in \calK,$ in $\calK: \ 
\calK=\bigsqcup_{g'\in \calG'}g'\calG_0.$  
Then we have $\calK=\calG'\calG_0,$ and 
$\calG'$ is a complete set of representatives of the double cosets 
$\calH g\calG_0, \ g \in \calG,$ in $\calG$ 
as well (Remark \ref{rem:twoG's}). 
Thus $\calK$ can serve as $\calL.$ 
For instance, consider the example immediately after Corollary
\ref{co:KZ}---the actions related to the two-sample Wishart problem.
Then $O(p)$ can serve as $\calL.$

\

Now we have by Theorem \ref{th:form of V} and \eqref{eq:calL condition} 
that
\[
\calV = \{ g_i \calG_0: i \in I \}
=\{ l \calG_0: l \in \calL \}
\]
is a global cross section for action $(\calH, \calG/\calG_0),$
and thus we obtain the decomposition
\begin{eqnarray}
\label{X<->UxVxZ}
\calX & \leftrightarrow & \calU \times \calV \times \calZ, 
\qquad  \calU = \calH/\calH_0, \nonumber \\
x & \leftrightarrow & (u, v, z), \ \ \ x = h l z,
\ \ u = h \calH_0, \ v = l \calG_0,
\end{eqnarray}
where $\calH_0 = \calH \cap \calG_0,$ since $\calL$ 
contains $e \in \calG$ and thus $\calG_0\in \calV.$ 
When $\calG'$ is taken in such a way that $e\in\calG',$ we have 
$\calG_0 < \calG'\calG_0=\calL$ and thus 
$\calV=\calL/\calG_0$ is the left coset space. 
We can always take $\calG'$ in this way, 
and we decide to do so. 

Concerning topological questions, we make the following 
assumptions in addition to (a) through (d) at the beginning of 
Section \ref{sec:decomposable-distribution}.

\begin{assumption} \quad 
\begin{enumerate}
\setlength{\itemsep}{0pt} 
 \item $\calH$ and $\calL$ are closed subgroups of $\calG.$
 \item $\calF$ is compact.
\end{enumerate}
\end{assumption}
Note that under our assumptions, $\calH_0 = \calH \cap \calG_0$
is compact, since $\calG_0$ is compact and $\calH \cap \calG_0$ is
closed in the relative topology of $\calG_0.$
Note also that the one-to-one correspondence
$\calG \leftrightarrow (\calH \times \calL)/\calF^*$
is a homeomorphism since $\calG$ is second 
countable (p.92 of Wijsman (1990)). 

As before, let $\lambda$ be a relatively invariant 
measure on $\calX$ under the action of $\calG$ 
with multiplier $\chi.$
Now we define the {\it extended decomposable distributions} as follows.

\begin{definition} \quad
A distribution on $\calX$ is said to be an extended 
decomposable distribution with respect to 
a pair of global cross sections $(\calZ, \calV)$ 
if it is of the form
\[
  f(x)\lambda(dx) =
  f_{ \calU }(u(x)) f_{ \calV }(v(x)) f_{ \calZ }(z(x)) \lambda(dx).
\]
\end{definition}

The following theorem gives the distributions of $u, \ v$ and $z$ when $x$
is distributed according to an extended decomposable distribution.

\begin{theorem} \quad
\label{th:distributions of u,v,z}
Suppose that $x$ is distributed according to an extended 
decomposable distribution
$f_{ \calU }(u(x)) f_{ \calV }(v(x)) f_{ \calZ }(z(x)) \lambda(dx).$
Then $u = u(x) = h \calH_0, \ v = v(x) = l \calG_0$ and
$z = z(x) \ (x = h l z)$ are independently distributed
with the joint distribution
\begin{eqnarray*}
& & f_{ \calU }(u) \chi(u) \mu_{ \calU }(du) \\
& & \ \ \ \times f_{ \calV }(v) \chi(v)
    \Delta^{ \calG }(v) \Delta^{ \calL}(v)^{-1}
    \mu_{ \calV }(dv) \\
& & \ \ \ \times f_{ \calZ }(z) \nu_{ \calZ }(dz),
\end{eqnarray*}
where $\Delta^{\calG}$ (resp. $\Delta^{\calL}$\kmms) is the right-hand modulus of
$\calG$ (resp. $\calL$\kmms),
measure $\mu_{ \calU }$ (resp. $\mu_{ \calV }$\kmms) is a version of the
invariant measures on $\calU = \calH/\calH_0$
(resp. $\calV = \calL/\calG_0$\kmms),
and $\nu_{ \calZ }$ is the probability measure 
in \eqref{eq:decomposition of lambda}.
\end{theorem}
This theorem can be proved by Proposition 7.6.1 and (7.6.5) of 
Wijsman (1990).

\section{Examples}

In this section, we apply Theorem \ref{th:distributions of u,v,z}
to the star-shaped distributions, the two-sample Wishart 
problem and the distributions of rankings. 

\subsection{Star-shaped distributions with symmetry}
Consider the actions related to the star-shaped distributions---the
actions of $\calG = \bbR^*_\times$ and $\calH = \bbR_+^*$ on
$\calX = \bbR^p - \{ {\bf 0} \}.$

Then we have
$\calG_0 = \calH_0 = \{ 1 \}.$
If we take
$\calL = \calG' 
= \{ \pm 1 \},$
we obtain the bijection
\begin{eqnarray*}
\calX & \leftrightarrow &
           \calH \times \calL \times \calZ, \\
\vx & \leftrightarrow & (h, \ \epsilon, \ \vz),
\ \ \ \vx = \epsilon h \vz,
\end{eqnarray*}
where $\calZ$ is a cross section for action $(\calG, \calX).$
Furthermore, $\calG=\bbR_{\times}^*$ and $\calL=\{ \pm 1 \}$ are unimodular:
$\Delta^{ \calG } = 1, \ \Delta^{ \calL } = 1.$

First, suppose that $\vx$ is distributed according to a star-shaped
distribution with respect to 
decomposable 
cross section $\tilde{\calZ}=\calG'\calZ 
=\calL \calZ
=\{ \epsilon \vz: \epsilon=\pm 1, \ \vz \in \calZ \}=\calZ \cup (-\calZ):$
\[
f(h(\vx))d\vx.
\]
This distribution can 
be regarded as the extended decomposable
distribution with
\begin{eqnarray*}
f_{ \calU }(h) &=& f_{ \calH }(h) = f(h), \\
f_{ \calV }( \epsilon ) &=& f_{ \calL }( \epsilon ) \equiv 1, \\
f_{ \calZ }( \vz ) & \equiv & 1,
\end{eqnarray*}
and
\[
\lambda(d\vx) = d\vx.
\]
Dominating measure $\lambda$ is relatively invariant under the action 
of $\calG$ with multiplier
\[
\chi(g) = |g|^p, \ \ \ g \in \calG,
\]
where $| \cdot |$ denotes the absolute value.
Therefore, we have by Theorem \ref{th:distributions of u,v,z} that
$h, \ \epsilon$ and $\vz$ are independently distributed according to
\begin{equation}
\label{eq:3parts} 
  \frac{1}{c_0} f(h) h^p h^{-1} dh 
  = \frac{1}{c_0} f(h) h^{p-1} dh, \quad
  \mu_{ \calV } = \mu_{ \calL } \quad {\rm and} \quad 
  \nu_{ \calZ },
\end{equation}
respectively, where $c_0 = \int_0^{ \infty } f(h) h^{p-1} dh$
and
$\mu_{ \calL }( \{ 1 \} ) = \mu_{ \calL }( \{ -1 \} )
= 1/2.$ 
We can see that the distributions of $\vx$ and $-\vx$ are the same 
by $-\vx \leftrightarrow (h, -\epsilon, \vz)$ for 
$\vx \leftrightarrow (h, \epsilon, \vz).$ 
Under the additional assumption that $h( \vx )$ is 
piecewise of class $C^1,$
we have
\[
\nu_{ \calZ }( d\vz ) = 2c_0 \langle \vz, \vn_{ \vz } \rangle d\vz,
\]
where $\vn_{ \vz }$ is the outward unit normal vector of $\calZ,$
$d\vz$ on the right-hand side is the volume element of $\calZ$
and $\langle \ \cdot \ , \ \cdot \ \rangle$ denotes the standard inner 
product (Section 4 of Kamiya, Takemura and Kuriki (2006)).

Next, by taking 
nonconstant $f_{\calL}(\epsilon)$ instead,  
we can make an asymmetric distribution of $\vx$ as follows. 
Suppose 
\begin{equation}
\label{eq:cfdx}
\vx\sim c(\epsilon(\vx))f(h(\vx))d\vx, 
\end{equation}
where 
\begin{equation}
\label{eq:c(e)}
c(\epsilon):=\begin{cases}
				c & \text{if $\epsilon =1,$} \\ 
				2-c & \text{if $\epsilon = -1$} 	
			\end{cases}
\end{equation}  
with $0\le c \le 2.$ 
Then again Theorem \ref{th:distributions of u,v,z} implies 
that $h, \ \epsilon$ and $\vz$ are independently distributed as 
\eqref{eq:3parts}, but this time with 
$\mu_{\calV}=\mu_{\calL}$ replaced by  
$\tilde{\mu}_{\calL}(\{1\})=c/2, \ 
\tilde{\mu}_{\calL}(\{-1\})=1-(c/2).$ 
In this case, the distribution of $-\vx$ is different from that of 
$\vx$ unless $c=1.$ 

Note that we can make a more general skewed 
distribution 
by considering a cross-sectionally contoured distribution 
$\tilde{f}(\epsilon(\vx)h(\vx))d\vx$ 
with respect to $\calZ.$  
But in that case, we have to specify the function $\tilde{f}$ defined 
on the whole of $\bbR_{\times}=
(-\infty, 0)\cup(0, \infty).$ 
By contrast, in the case of \eqref{eq:cfdx} we have only to 
specify $f$ defined on $\bbR_{+}=(0, \infty)$ and one value $c$ in 
\eqref{eq:c(e)}. 


\subsection{Two-sample Wishart problem}

Consider the action of $\calG=GL(p)$ on the sample space 
$\calX$ in \eqref{eq:wishart sample space}. 
The cross section $\calZ$ is taken as \eqref{eq:wishart cross section} with 
the isotropy subgroup $\calG_0$ in \eqref{eq:G0 for wishart}.

We continue to take $\calH=LT(p)$ and 
$\calL = \calK = O(p).$
Then we obtain the bijection
\begin{eqnarray*}
\calX & \leftrightarrow &
           \calH \times \calL/\calG_0 \times \calZ, \\
(W_1, W_2) & \leftrightarrow & 
            \left(T, \ C\calG_0, \ (\Lambda, \ I_p - \Lambda)\right),
\ \ \ (W_1, W_2) =
\left(T C \Lambda C^\trans T^\trans, \ 
T C (I_p - \Lambda) C^{ \trans} T^\trans \right).
\end{eqnarray*}
Furthermore, $\calG=GL(p)$ and $\calL=O(p)$ are unimodular:
$\Delta^{ \calG } = 1, \ \Delta^{ \calL } = 1.$

Suppose that the random matrices $W_1$ and $W_2$ are independently
distributed according to $W_p(n_1, \Sigma)$ and $W_p(n_2, \Sigma),$
respectively.
Then the distribution of $(W_1, W_2)$ can be regarded as the extended
decomposable distribution with
\[
f_{ \calU }(T) = f_{\calH}(T) 
\propto {\rm \ etr}\left(-\frac{1}{2}\Sigma^{-1}TT^\trans \right),
\ \ \ 
f_{ \calV }(C \calG_0 ) = f_{\calL/\calG_0}(C\calG_0) \equiv 1,
\ \ \ 
f_{ \calZ }\left((\Lambda, \ I_p - \Lambda)\right) \equiv 1
\]
and 
\begin{equation}
\label{eq:lambda for wishart}
\lambda(d(W_1,W_2))
=(\det W_1)^{a-\frac{p+1}{2}} (\det W_2)^{b-\frac{p+1}{2}}dW_1 dW_2, 
\end{equation}
where $a=n_1/2, \ b=n_2/2, \ W_1=(w_{1,ij}), \ W_2=(w_{2,ij}), \ 
dW_1=\prod_{i\ge j}dw_{1,ij}, \ dW_2=\prod_{i\ge j}dw_{2,ij}.$
The dominating measure $\lambda(d(W_1, W_2))$ is relatively 
invariant with multiplier
\[
\chi(B)=(\det B)^{2(a+b)}=(\det B)^{n_1+n_2}, \quad B \in GL(p), 
\]
(Wijsman (1990), (9.1.4)).
It follows from Theorem \ref{th:distributions of u,v,z} that    
$T, \ C \calG_0$ and $\Lambda$ are independently distributed.
The distributions of these parts are given in standard textbooks of
multivariate statistical theory (see Anderson (2003) or
Muirhead (1982), for example).
In particular, $C\calG_0$ is distributed according to the
invariant probability measure on $\calV = \calL/\calG_0=O(p)/\calG_0$ 
induced by the 
Haar measure on $O(p).$

A nonstandard distribution is given as follows.
Since the normalizer of $\calG_0$ is 
\begin{align*}
\calN = \{ P \in GL(p) : 
& \ P \ \text{has exactly one nonzero element} \\
& \quad \quad \text{in each row and in each column} \}, 
\end{align*}
we know from \eqref{eq:g0nz} that a general global cross section 
$\calZ'$ is of the form
\begin{align}
\label{eq:arbitrary cross section for two-sample Wishart}
\calZ'=
&  \Big\{\left(B P(\Lambda) \Lambda P(\Lambda)^\trans B^\trans, 
\ B P(\Lambda)(I_p-\Lambda) P(\Lambda)^\trans B^\trans \right): \\
& \qquad\qquad
  \ \ \Lambda = \diag (\lambda_1,\ldots,\lambda_p),
  \ 1 > \lambda_1 > \cdots > \lambda_p > 0 \Big\} \nonumber
\end{align}
with $B \in GL(p)$ and $P( \Lambda ) \in \calN.$ 
Without loss of generality, we assume $B=I_p$ in 
\eqref{eq:arbitrary cross section for two-sample Wishart}. 
Let $B(W)$ denote the equivariant part of $W=(W_1, W_2)$ with respect to 
$\calZ.$  
Then the equivariant part with respect to $\calZ'$ is 
$B(W)P(\Lambda(W))^{-1}\calG_0$ by \eqref{eq:transformation of y}.
Writing the latter as 
\[
  B(W)P(\Lambda(W))^{-1}\calG_0 = T'(W) C'(W)\calG_0, \qquad 
  T'(W) \in LT(p), \ C'(W) \in O(p), 
\]
we obtain the decomposition 
\[
W \leftrightarrow \left( T'(W), \ C'(W)\calG_0, \ 
z'(W) \right), 
\]
where  
$z'(W)=(P(\Lambda(W)) \Lambda(W) P(\Lambda(W))^\trans, \ 
P(\Lambda(W))(I_p-\Lambda(W)) P(\Lambda(W))^\trans )\in \calZ'$ 
is the invariant part of $W$ with respect to $\calZ'.$ 
Suppose that the density $f(W)$ with respect 
to $\lambda$ in \eqref{eq:lambda for wishart} 
with general $a, b>(p+1)/2$ is 
factored as 
\[
f(W)=f_{\calH}(T'(W))f_{\calL/\calG_0}(C'(W)\calG_0)f_{\calZ'}(z'(W)).
\]  
Then we can get from Theorem \ref{th:distributions of u,v,z}
the distributions of the three parts 
$T'=T'(W), \ C'\calG_0=C'(W)\calG_0$ and $z'=z'(W)$ 
using 
$
\Delta^{\calG}=\Delta^{\calL}=1, \ 
\chi(B)=(\det B)^{2(a+b)}$ 
and 
$\mu_{\calH}(dT)=
\mu_{LT(p)}(dT) 
= 
\prod_{i=1}^p t_{ii}^{-i}dT, \ dT=\prod_{i \ge j}dt_{ij}, \ 
T=(t_{ij})$ 
(Wijsman (1990), (7.7.2)).  
Specifically, 
\begin{eqnarray*}
T' &\sim& \frac{1}{c_{\calH}}f_{\calH}(T')\prod_{i=1}^p t_{ii}'^{2a+2b-i}dT', \\ 
C'\calG_0 &\sim& f_{\calL/\calG_0}(C'\calG_0)d\mu_{\calL/\calG_0}(C'\calG_0), \\ 
z' &\sim& \frac{1}{c_{\calZ'}}f_{\calZ'}(z')d\nu_{\calZ'}(z'), 
\end{eqnarray*}
where $c_{\calH}=\int_{LT(p)}f_{\calH}(T)\prod_{i=1}^p t_{ii}^{2a+2b-i}dT, \ 
c_{\calZ'}=\int_{\calZ'}f_{\calZ'}(z)d\nu_{\calZ'}(z)$ 
and $\mu_{\calL/\calG_0}$ is an appropriate version of the invariant measures on 
$\calL/\calG_0=O(p)/\calG_0.$

\subsection{Decompositions of rankings}

Our general discussion can be applied to discrete distributions as well. 
In this subsection, let us look at one such example---distributions of 
preference rankings. 
For the analysis of ranking data in general, the reader is referred to the excellent books by 
Critchlow (1985), Diaconis (1988) and Marden (1995). 
Other interesting problems about preference rankings can be found 
in Kamiya, Orlik, Takemura and Terao (2006). 

Let us consider rankings of $m$ objects $1, 2, \ldots, m.$
We denote rankings as $\sigma=(\sigma(1), \sigma(2), \ldots, \sigma(m)),$ 
where $\sigma(i)$ stands for the rank given to object $i.$
Then we can regard 
\[
\sigma=(\sigma(1), \sigma(2), \ldots, \sigma(m))=
\begin{pmatrix}
1 & 2 & \cdots & m \\ 
\sigma(1) & \sigma(2) & \cdots & \sigma(m)
\end{pmatrix}
\] 
as an element of the symmetric group $S_m$ on 
$\{ 1, 2, \ldots, m \}.$ 
We can deal with distributions of rankings $\sigma \in S_m$ 
by considering 
probability functions on 
$\calX=S_m.$ 

Here we define an action of $\calG=S_{m-1}<S_m$ 
on $\calX=S_m$ as follows. 
Thinking of $\tau\in S_{m-1}$ as permutations of 
$\{ 2, \ldots, m \},$ 
we express $\tau$ as 
\[
\tau=(\tau(2), \ldots, \tau(m))=
\begin{pmatrix}
2 & \cdots & m \\ 
\tau(2) & \cdots & \tau(m)
\end{pmatrix}.
\] 
Permutation $\tau\in S_{m-1}$ changes rank 
$k\in \{2,\ldots,m\}$ to rank $\tau(k)\in \{2,\ldots,m\}.$ 
Now we consider the action $(\tau, \sigma) \mapsto \tau\sigma.$
Here $\tau\sigma\in S_m$ means 
$(\tau\sigma)(i)=\tau(\sigma(i)), \ i\in \{1,2,\ldots, m\},$ 
where we agree that $\tau(1)=1.$ 
Note that this action is free:$\calG_0=\{ e \}, \ e=(2,\ldots,m)
=(1,2,\ldots,m).$ 

Under the above action, the orbit containing $\sigma\in S_m$ is 
\[
\calG\sigma=S_{m-1}\sigma=\{ \tilde{\sigma}\in S_m: 
\tilde{\sigma}^{-1}(1)=\sigma^{-1}(1)\}
\] 
and $\calX=S_m$ consists of the $m$ orbits 
\[
\{\tilde{\sigma} \in S_m: \tilde{\sigma}(i)=1 \}, \quad  
i=1,2,\ldots,m.
\] 
Orbit $\{\tilde{\sigma} \in S_m: \tilde{\sigma}(i)=1 \}$ is seen to be  
the set of 
all rankings that rank object $i$ first, and we will call this the 
$i$\kmms-orbit.

One way of selecting a representative $\sigma_i$ 
of the $i$\kmms-orbit is choosing $\sigma_i$ as 
\begin{equation}
\label{eq:2<cdots<m}
\sigma_i^{-1}(1)=i, \quad \sigma_i^{-1}(2)<\cdots<\sigma_i^{-1}(m).
\end{equation}
Let us take the cross section $\calZ$ consisting of these representatives 
$\sigma_i, \ i=1,2,\ldots,m,$ i.e.,   
$\calZ=\{ \sigma_i: i=1,2,\ldots, m \}.$ 
Then, with respect to this cross section, we can write an arbitrary 
$\sigma\in \calX=S_m$ uniquely as 
\[
\sigma=\tau s=\tau(\sigma) s(\sigma), \qquad \tau\in \calG=S_{m-1}, \ \  
s\in \calZ.
\] 
Note that $s(\sigma)$ can be expressed explicitly as 
$s(\sigma)=\sigma_{\sigma^{-1}(1)}.$ 
For example, when $m=4,$ ranking   
\[
\sigma=(4,2,1,3)=
\begin{pmatrix}
1 & 2 & 3 & 4 \\ 
4 & 2 & 1 & 3
\end{pmatrix}
\] 
belongs to the $3$\kmms-orbit, and the corresponding 
$\tau(\sigma)$ and $s(\sigma)$ are given as 
\[
\tau(\sigma)=
(4,2,3)=
\begin{pmatrix}
2 & 3 & 4 \\ 
4 & 2 & 3
\end{pmatrix}, \quad 
s(\sigma)=\sigma_3=
(2,3,1,4)=
\begin{pmatrix}
1 & 2 & 3 & 4 \\ 
2 & 3 & 1 & 4
\end{pmatrix}.
\]

We are now in a position to introduce a family of 
distributions on $\calX=S_m.$ 
Consider the following type of probability functions:
\begin{equation}
\label{eq:p=pp}
p(\sigma)=\tilde{p}_{\calG}(d(\sigma, s(\sigma)))p_{\calZ}(s(\sigma)), \quad 
\sigma\in S_m,
\end{equation}
where $d(\ \cdot \ , \ \cdot \ )$ is a right-invariant metric on 
$S_{m-1}$ 
(note that $\sigma^{-1}(1)=s(\sigma)^{-1}(1)$\kmms).  
We usually take $\tilde{p}_{\calG}$ to be a decreasing function. 
In that case, 
$\tilde{p}_{\calG}(d(\sigma, s(\sigma)))$ expresses a unimodal 
distribution in the orbit with the representative 
$s(\sigma)$ as the modal ranking; 
the farther away from the representative $s(\sigma),$ 
the smaller probability of that ranking $\sigma.$ 
An example of a specification of $\tilde{p}_{\calG}$ is provided by 
Mallows' model $\tilde{p}_{\calG}(d(\sigma, s(\sigma)))=
c\exp(\theta d(\sigma, s(\sigma)))$ with $\theta \le 0,$ where $c$ is 
the normalizing constant (Mallows (1957), Feigin and Cohen (1978), 
Fligner and Verducci (1986)).  
Now, writing $p_{\calG}( \tau ):=\tilde{p}_{\calG}(d(\tau, e)),$ 
we can express \eqref{eq:p=pp} as  
\[
p(\sigma)=p_{\calG}(\tau(\sigma))p_{\calZ}(s(\sigma)), \quad 
\sigma\in S_m.
\]  
This can be regarded as a decomposable distribution  
under our action of $\calG=S_{m-1}$ on $\calX=S_m.$ 
We are writing $p(\sigma)$ instead of $f(\sigma)$ in order 
to emphasize that we are dealing with discrete probability functions. 
Of course, the density is with respect to the counting measure, 
so $\chi=1.$ 

In the discussions so far, the representatives 
$\sigma_1,\ldots,\sigma_m$  
of the orbits have been chosen as \eqref{eq:2<cdots<m},
but this is just one way of choosing them and other 
choices are also possible. 
The choice is arbitrary as long as they form a complete 
set of representatives.  
In fact, the selection rule can differ from orbit to orbit. 
For example, we may consider taking the representative 
$\sigma_{i_0} \in \{\tilde{\sigma} \in S_m: \tilde{\sigma}(i_0)=1 \}$ 
for some object $i_0 \in \{ 1,2,\ldots,m \}$ as 
follows: 
for some object $j_0,$ 
\begin{equation}
\label{eq:sigma_i0}
\sigma_{i_0}^{-1}(1)=i_0, 
\quad \sigma_{i_0}^{-1}(2)<\cdots<\sigma_{i_0}^{-1}(m-1), 
\quad \sigma_{i_0}^{-1}(m)=j_0.
\end{equation}
For the other orbits, we continue to select 
representatives $\sigma_i, \ i \ne i_0,$ as in \eqref{eq:2<cdots<m}. 


We now discuss motivations of the above modeling. 
Imagine we are considering people's preference rankings of 
a league of $m$ sports teams $1,2,\ldots,m.$
Suppose that we are faced with the following situation: 
while people are interested in their favorite team,   
they do not care much about the differences among the rest 
and tend to simply rank the second to last preferred  
teams according to the ranks in the standings (based on 
winning percentages).
In that case, we can describe the situation by 
modeling the distribution of the rankings as follows. 

We label the $m$ teams 
with the ranks in the standings.
Then, we choose the representatives of all the orbits as in  
\eqref{eq:2<cdots<m}.
With these choices, the distribution \eqref{eq:p=pp} of people's 
preference rankings 
implies that their top rank is distributed 
according to $p_{\calZ}(s(\sigma)),$ while the rest of the ranks 
are distributed 
based on the distances to the 
modal ranking, 
which in this case is 
the ranking (having the same relative ranks of 
the non-top objects as the ranking) in the standings; 
the closer to the ranking in the standings, the larger percentage of 
people with that ranking. 
So we have succeeded in describing the case in question.

Alternatively, there may be some cases where 
the fans of a certain team $i_0$ have a strong sense of rivalry with 
some other team $j_0.$
In those cases, 
the choice \eqref{eq:sigma_i0} will be more appropriate. 
%


Now, we can go further and consider decompositions into 
three parts. 
Suppose that people are very interested in the top rank,    
interested to some degree in 
ranks $2,\ldots,m' \ 
(2 \le m' \le m-1)$  
and totally indifferent about 
the rest of the ranks $m'+1,\ldots,m.$ 

Let us consider the action of 
$\calH=S_{m-m'}$ on $
\calG=S_{m-1},$ 
with $S_{m-m'}$ 
being regarded as the set of permutations of 
$\{ m'+1,\ldots,m \}.$ 
Denote a cross section of this action by $\calV.$ 
Then any $\tau\in \calG$ can be written uniquely as 
$\tau=ht, \ h \in \calH, \ t \in \calV.$ 
Hereafter we will write $
S_{m}=S_{\{1,\ldots,m\}}, \ 
S_{m-1}=S_{\{2,\ldots,m\}}$ and 
$
S_{m-m'}=S_{\{m'+1,\ldots,m\}}.$

By writing $\sigma \in \calX=S_{\{ 1,\ldots,m\}}$ as 
\begin{eqnarray*}
\sigma
&=& \tau s \quad 
(\tau \in \calG=S_{\{2,\ldots,m\}}, \ s \in \calZ) \\ 
&=& h t s \quad 
(h \in \calH=S_{\{m'+1,\ldots,m\}}, \ t \in \calV), 
\end{eqnarray*}
we obtain bijections 
\begin{eqnarray}
\label{eq:bijections}
\sigma 
& \leftrightarrow &  
\left( \tau, \ S_{\{2,\ldots,m\}}\sigma \right) 
 \ (\in S_{\{2,\ldots,m\}}\times 
S_{\{2,\ldots,m\}}\backslash S_{\{1,\ldots,m\}}) \nonumber \\ 
& \leftrightarrow &  
\left(
h, \ S_{\{m'+1,\ldots,m\}}\tau, \ S_{\{2,\ldots,m\}}\sigma 
\right) \nonumber \\ 
&& \qquad (\in S_{\{m'+1,\ldots,m\}} \times 
S_{\{m'+1,\ldots,m\}}\backslash S_{\{2,\ldots,m\}}\times 
S_{\{2,\ldots,m\}}\backslash S_{\{1,\ldots,m\}}) \\ 
& \leftrightarrow & 
\left(h, \ t, \ s\right)
\in S_{\{m'+1,\ldots,m\}}\times \calV \times \calZ. \nonumber 
\end{eqnarray} 
Remember that $\calH \backslash \calG$ 
for $\calH<\calG$ denotes the right coset space: 
$\calH \backslash \calG=\{ \calH g: g \in \calG \}.$  

\begin{example} \quad 
Let $m=6, \ m'=3,$ and take 
$\calZ=\{ \sigma_i: i =1,2,\ldots, 6 \}$ with 
$\sigma_1,\ldots,\sigma_6$ as in \eqref{eq:2<cdots<m} and 
$\calV=\{ \sigma_{(i, j)}\in S_{\{ 2,\ldots,6 \}}: 
\{ i, j \} \subset \{ 2,\ldots,6 \}, \ i\ne j \}$ with 
$\sigma_{(i, j)}$ 
such that $\sigma_{(i, j)}^{-1}(2)=i, \ 
\sigma_{(i, j)}^{-1}(3)=j, \ 
\sigma_{(i, j)}^{-1}(4)<\sigma_{(i, j)}^{-1}(5)<\sigma_{(i, j)}^{-1}(6).$  
Then, 
\begin{eqnarray*}
\begin{pmatrix}
1 & 2 & 3 & 4 & 5 & 6 \\ 
6 & 3 & 5 & 1 & 4 & 2 
\end{pmatrix} 
&=& 
\begin{pmatrix}
2 & 3 & 4 & 5 & 6 \\ 
6 & 3 & 5 & 4 & 2  
\end{pmatrix}
\begin{pmatrix}
1 & 2 & 3 & 4 & 5 & 6 \\ 
2 & 3 & 4 & 1 & 5 & 6 
\end{pmatrix} \\ 
&=& 
\begin{pmatrix}
4 & 5 & 6 \\ 
6 & 5 & 4   
\end{pmatrix} 
\begin{pmatrix}
2 & 3 & 4 & 5 & 6 \\ 
4 & 3 & 5 & 6 & 2  
\end{pmatrix}
\begin{pmatrix}
1 & 2 & 3 & 4 & 5 & 6 \\ 
2 & 3 & 4 & 1 & 5 & 6 
\end{pmatrix},
\end{eqnarray*}
so
\begin{eqnarray*}
&& 
S_{\{1,2,3,4,5,6\}}
\ni 
\begin{pmatrix}
1 & 2 & 3 & 4 & 5 & 6 \\ 
6 & 3 & 5 & 1 & 4 & 2 
\end{pmatrix} \\ 
& \leftrightarrow & 
\left( 
\begin{pmatrix}
2 & 3 & 4 & 5 & 6 \\ 
6 & 3 & 5 & 4 & 2  
\end{pmatrix}, \ 
S_{\{2,3,4,5,6\}}
\begin{pmatrix}
1 & 2 & 3 & 4 & 5 & 6 \\ 
6 & 3 & 5 & 1 & 4 & 2 
\end{pmatrix}
\right) \\ 
&& \qquad \qquad 
(\in S_{\{ 2,3,4,5,6\}}\times S_{\{2,3,4,5,6\}}\backslash S_{\{1,2,3,4,5,6\}})
\\ 
& \leftrightarrow & 
\left(\
\begin{pmatrix}
4 & 5 & 6 \\ 
6 & 5 & 4   
\end{pmatrix}, \ 
S_{\{4,5,6\}} 
\begin{pmatrix}
2 & 3 & 4 & 5 & 6 \\ 
6 & 3 & 5 & 4 & 2  
\end{pmatrix}, \ 
S_{\{2,3,4,5,6\}}
\begin{pmatrix}
1 & 2 & 3 & 4 & 5 & 6 \\ 
6 & 3 & 5 & 1 & 4 & 2 
\end{pmatrix}
\right) \\ 
&& \qquad \qquad \qquad 
(\in S_{\{ 4,5,6 \}} \times 
S_{\{ 4,5,6 \}}\backslash S_{\{2,3,4,5,6 \}} \times 
S_{\{ 2,3,4,5,6 \}}\backslash S_{\{ 1,2,3,4,5,6 \}}) \\ 
& \leftrightarrow & 
\left(
\begin{pmatrix}
4 & 5 & 6 \\ 
6 & 5 & 4   
\end{pmatrix}, \ 
\begin{pmatrix}
2 & 3 & 4 & 5 & 6 \\ 
4 & 3 & 5 & 6 & 2  
\end{pmatrix}, \ 
\begin{pmatrix}
1 & 2 & 3 & 4 & 5 & 6 \\ 
2 & 3 & 4 & 1 & 5 & 6 
\end{pmatrix}
\right) 
\in S_{\{4,5,6\}}\times \calV \times \calZ.
\end{eqnarray*}
\end{example}

Now consider probability functions $p(\sigma)$ which can be factored  
with respect to 
\eqref{eq:bijections} as  
\begin{equation}
\label{eq:p=pVpZ}
p(\sigma)=
p_{\calH\backslash\calG}(S_{\{m'+1,\ldots,m\}}\tau) 
p_{\calZ}(s)
\end{equation}
with 
$p_{\calH\backslash\calG}(S_{\{m'+1,\ldots,m\}} \tau)
:=\tilde{p}_{\calH\backslash\calG}
(d'(S_{\{m'+1,\ldots,m\}}\tau, S_{\{m'+1,\ldots,m\}})),$ 
where $d'$ is the 
Hausdorff metric on 
$\calH\backslash\calG=S_{\{m'+1,\ldots,m\}}\backslash S_{\{2,\ldots,m\}}$ 
induced by the metric $d$ on $\calG=S_{\{2,\ldots,m\}}$ 
(Critchlow (1985)). 
Note that distribution \eqref{eq:p=pVpZ} can be seen as an 
extended decomposable distribution 
with respect to 
\eqref{eq:bijections}. 
Now, label the teams 
with the ranks in the standings and 
take $\calZ=\{ \sigma_i: i=1,2,\ldots, m \}$ with $\sigma_i, \ i=1,2,\ldots,m,$ in  
\eqref{eq:2<cdots<m} as before.  
Then \eqref{eq:p=pVpZ} describes the situation stated earlier. 
Note that in \eqref{eq:p=pVpZ} the choice of $\calV$ is irrelevant,
because the rest of the ranks are uniform.

\begin{rem}
Here we have studied the permutations of ranks. 
By considering the orderings 
\[
(\text{top ranked object, second ranked object, \ldots})
\] 
instead of rankings 
\[
(\text{object 1's rank, object 2's rank, \ldots}), 
\] 
we can also deal with the permutations of objects in a similar manner. 
\end{rem}


\begin{thebibliography}{99}

\bibitem{Anderson2003} 
Anderson, T. W. (2003). 
{\it An Introduction to Multivariate Statistical Analysis}, Third Edition. 
John Wiley \& Sons, New Jersey.

\bibitem{Barndorff-Nielsen--Blaesild--Eriksen} 
Barndorff-Nielsen, O. E., Bl{\ae}sild, P.\ and Eriksen, P. S. (1989). 
{\it Decomposition and Invariance of Measures, and Statistical 
Transformation Models}. 
Lecture Notes in Statistics, Vol.\ 58, 
Springer-Verlag, Berlin. 

\bibitem{Bondar}
Bondar, J. V. (1976). 
Borel cross-sections and maximal invariants. 
{\it Ann. Statist.}, {\bf 4}, 866--877.

\bibitem{Cambanis-Huang-Simons}
Cambanis, S., Huang, S.\ and Simons, G. (1981). 
On the theory of elliptically contoured distributions. 
{\it J. Multivariate Anal.}, {\bf 11}, 368--385.

\bibitem{Critchlow1985} 
Critchlow, D. E. (1985). 
{\it Metric Methods for Analyzing Partially Ranked Data}. 
Lecture Notes in Statistics, Vol.\ 34, 
Springer-Verlag, Berlin. 

\bibitem{Diaconis} 
Diaconis, P. (1988). 
{\it Group Representations in Probability and Statistics}.  
Lecture Notes--Monograph Series, Vol.\ 11, 
Institute of Mathematical Statistics, Hayward, California.

\bibitem{Eaton} 
Eaton, M. L. (1989).  
{\it Group Invariance Applications in Statistics}.  
Regional Conference Series in Probability and Statistics, Vol.\ 1, 
Institute of Mathematical Statistics, Hayward, California.

\bibitem{Fang-Anderson}
Fang, K. T.\ and Anderson, T. W. (1990).  
{\it Statistical Inference in Elliptically Contoured and Related Distributions}. 
Allerton Press, New York.

\bibitem{Fang-Zhang} 
Fang, K. T.\ and Zhang, Y. T. (1990). 
{\it Generalized Multivariate Analysis}. 
Springer-Verlag, New York.

\bibitem{Farrell} 
Farrell, R. H. (1985).  
{\it Multivariate Calculation: Use of the Continuous Groups}. 
Springer-Verlag, New York. 

\bibitem{Feigin-Cohen} 
Feigin, P. D.\ and Cohen, A. (1978). 
On a model for concordance between judges. 
{\it J. Roy. Statist. Soc.} B, {\bf 40}, 203--213.

\bibitem{Fernandez-Osiewalski-Steel}
Fern\'{a}ndez, C., Osiewalski, J.\ and Steel, M. F. J. (1995). 
Modeling and inference with $v$\kmms-spherical distributions.
{\it J. Amer. Statist. Assoc.}, {\bf 90}, 1331--1340.

\bibitem{Ferreira-Steel}
Ferreira, J. T. A. S.\ and Steel, M. F. J. (2005). 
Modelling directional dispersion through hyperspherical log-splines. 
{\it J. Roy. Statist. Soc.} B, {\bf 67}, 599--616.

\bibitem{Fligner-Verducci} 
Fligner, M. A.\ and Verducci, J. S. (1986). 
Distance based ranking models. 
{\it J. Roy. Statist. Soc.} B, {\bf 48}, 359--369. 

\bibitem{Kamiya}
Kamiya, H. (1996). 
Borel isomorphism between the sample space and a 
product space. 
{\it Math. Methods Statist.}, {\bf 5}, 237--243.

\bibitem{Kamiya-Orlik-Takemura-Terao} 
Kamiya, H., Orlik, P., Takemura, A.\ and Terao, H. (2006). 
Arrangements and ranking patterns. 
{\it Ann. Comb.}, {\bf 10}, 267-284. 

\bibitem{Kamiya-Takemura-Kuriki}
Kamiya, H., Takemura, A.\ and Kuriki, S. (2006).
Star-shaped distributions and their generalizations. 
{\it J. Statist. Plann. Inference}, to appear. 
Available at {\bf http://arxiv.org/abs/math.ST/0605600}.

\bibitem{Kelker} 
Kelker, D. (1970).
Distribution theory of spherical distributions and a location-scale 
parameter generalization. 
{\it Sankhy\={a}} A, {\bf 32}, 419--430.

\bibitem{Koehn}
Koehn, U. (1970). 
Global cross sections and the densities of maximal invariants. 
{\it Ann. Math. Statist.}, {\bf 41}, 2045--2056.

\bibitem{Mallows}
Mallows, C. L. (1957). 
Non-null ranking models. I. 
{\it Biometrika}, {\bf 44}, 114--130. 

\bibitem{Marden} 
Marden, J. I. (1995). 
{\it Analyzing and Modeling Rank Data}. 
Chapman \& Hall, London.  

\bibitem{Muirhead} 
Muirhead, R. W. (1982).
{\it Aspects of Multivariate Statistical Theory}.
John Wiley \& Sons, New York. 

\bibitem{Osiewalski-Steel}
Osiewalski, J.\ and Steel, M. F. J. (1993). 
Robust Bayesian inference in $l_q$\kmms-spherical models. 
{\it Biometrika}, {\bf 80}, 456--460.

\bibitem{Wijsman1967}
Wijsman, R. A. (1967). 
Cross-sections of orbits and their application to 
densities of maximal invariants. 
{\it Proc. Fifth Berkeley Sympos. Math. Statist. and Probability}, 
Vol.~I: Statistics 389--400. 
Univ. California Press, Berkeley, California. 

\bibitem{Wijsman1986} 
Wijsman, R. A. (1986). 
Global cross sections as a tool for factorization of measures 
and distribution of maximal invariants. 
{\it Sankhy\={a}} A, {\bf 48}, 1--42.

\bibitem{Wijsman90} 
Wijsman, R. A. (1990).  
{\it Invariant Measures on Groups and Their Use in Statistics}.  
Lecture Notes--Monograph Series, Vol.\ 14, 
Institute of Mathematical Statistics, Hayward, California.

\end{thebibliography}
\end{document}